# Fractals in Africanist Music (and Dance)

Claudio Gómez-Gónzales, Sidhanth Raman, Siddharth Viswanath, Jesse Wolfson

The term "fractal" was introduced by Benoit Mandelbrot in 1975 [23] to describe geometric objects exhibiting self-similarity at scale. This concept, and the examples it encompasses, has given rise to rich theories, scientific and industrial applications, and works of art that make sense of previously unwieldy fragmentations occurring in natural phenomena. We often trace the lineage of these ideas to Weierstrass' continuous-but-non-differentiable curves and related pathological functions, yet we find diverse examples in abundance in Africanist material cultures (historically and contemporaneously), such as in textiles, architecture, and religious symbols [13]. Here, we consider how understanding, teaching, and doing mathematics engaged with African and African diasporic use of fractals provides an opening for embodied, experiential, and participatory engagement.

In 2012, choreographer Reggie Wilson was sourcing and developing material for a new concert dance work that would ultimately premier as *Moses(es)* in 2013. As Wilson described [36]:

> "Many of the movements and performance practices I was engaging with exhibited Africanist formal features--rhythms, brief sequences of movement, striking variations in quality and force or tempo--that did not map easily onto the schema and formalisms that characterize traditional Western concert dance. Around this same time, I encountered Ron Eglash's *African Fractals* which documents extensive use of fractal symmetry in African material cultures. I hypothesized that fractals were also prevalent in the performative cultures of Africa, and I invited Wolfson to help "translate" Eglash's text and to give his perspective on a selection of Africanist forms of music and dance that had most piqued my curiosity to possible usage of things fractal."

Wilson brought me (Jesse) in as a "math consultant" to assist the *Fist and Heel Performance Group* company members as they engaged with Africanist sources and developed original movement sequences for Wilson's well-received original concert dance *Moses(es)*. This experience confirmed the practical value of the hypothesis in the rehearsal studio: it provides tools for analyzing, playing with, reverse engineering, and generating movements in the service of creating original concert dance. The company developed a shared vocabulary for discussing and analyzing certain previously opaque features using the language of fractals, and Wilson and the dancers were able to reliably use this vocabulary to analyze and understand sourced movement in new and generative ways. The work described here grows directly out of Wilson's and the *Fist and Heel*'s "research into performance" practice, seeking to understand the mathematical underpinnings of Wilson's hypothesis through quantitative analysis of Africanist music samples.

## Fractals queer our sense of space

Dimension can be a tricky concept—we often speak of degrees of freedom, independent subsets, ascending chains, or even covering combinatorics—and fractals exceed our usual intuitions. While fractals are famous for evading a precise and consistent definition, we can characterize them in terms of their *Hausdorff dimension*, a quantity that gives us a sense of how "rough" a fractal might be. For those new to this concept, we can imagine dimension $d$ as teaching us how the "mass" $M$ of an object changes with respect to scaling by some factor $s$. Scaling by $s=2$ causes a solid cube to grow $8$ times heavier while a square quadruples its area and a line segment is simply twice as long. We write $M = s^d$ for this relationship, or $d = \log_s M$.

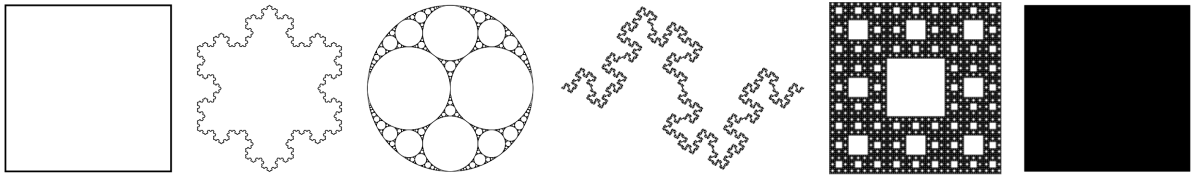

Figure 1: Examples of dimension. From left to right: the boundary of a square ($d = 1$, non-fractal), a Koch curve ($d = \log_3(4) = 1.2619...$), an Apollonian gasket ($d = 1.3057...$), a Minkowski sausage ($d = 1.5$), a Sierpinski carpet ($d = \log_3(8) = 1.8928...$), and a solid square ($d = 2$, non-fractal).

These paradoxical properties lead to many celebrated and exciting results [12], such as curve-like objects residing in a bounded area whose length can only be infinite. We encounter marvelous examples in which we might expect high dimensionality to be advantageous (for example, in our own lung alveoli [37]). Drawing from work of Richardson demonstrating that a coastline's length depends on the length of one's measuring stick, Mandelbrot [22] investigated precisely this phenomenon in applying fractional dimension to problems in nature. This sort of approach lends itself to methods of estimating fractal dimension.

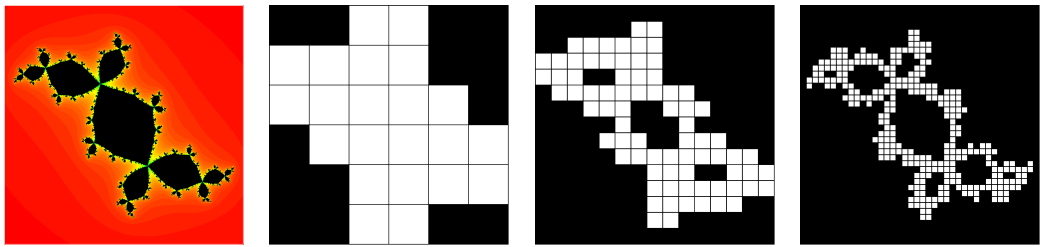

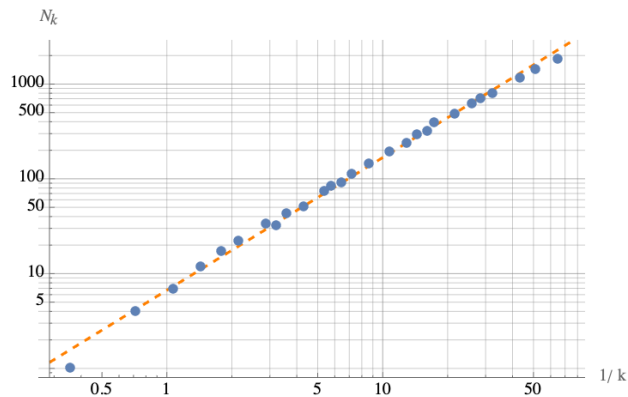

Figure 2: Computing the fractal dimension of a Douady rabbit (the leftmost image, $c = -0.123 + 0.745i$). Using boxes of size $k$ (from left to right, we have $k = 7/15, 14/75, 7/100$) we compare with the number needed to cover the rabbit's boundary. The dimension is approximated by plotting this relationship on log-log paper and reading off the best fit slope: here we measured $d \approx 1.3996$, very nearly the accepted $d = 1.3934...$ for this fractal.

A student learning logarithms—struck by the beauty, pervasiveness, and utility of fractals in the world around them—can realize the need for such a counterintuitive function, and similar wonder can inspire practitioners of analysis, dynamics or topology. Fractals also provide a context for us to question the conditions that give rise to mathematical activity and the natures of our practice itself. They exemplify the "lamentable scourge" [16] of counterexamples—Weierstrass' curves, Koch's snowflake, Cantor's dust, and more—dreamt up during the foundational mathematical crisis at the turn of the twentieth century. At the same time, they represent organizing principles in traditional knowledge systems of many Africanist cultures, such as the layout of Ba-ila kraals in southern Zambia, sculptural aesthetics from Lega societies of Zaire, and both urban planning and architectural design across Cameroon [13]. Fractals continue to be studied and implemented towards understanding space, assembling spiritual and communal life, solving practical problems, and creating art—as examples of participatory mathematical research extending beyond the template of professional mathematicians doing math to or on behalf of someone else.

## Hausdorff dimension and Africanist music

The present project began during the Spring of 2021 when we (Sidhanth, Siddharth, and Jesse) set out to explore Wilson's hypothesis on mathematical terms. This research involved an analysis of the music samples Wilson had originally shared, using our collective expertise in signal processing, sound engineering, and mathematics. We wondered if the passages that Wilson suspected to be fractal exhibited fractional Hausdorff dimensions (in particular, as waveforms analyzed as time-series data). We were also curious if listeners can reliably correlate our qualitative, acoustic perceptions of songs with the presence of fractional Hausdorff dimension.

To explore these questions, we developed software capable of computing Hausdorff dimension from audio files. Our approach uses Higuchi's method [17], which measures the signal "length" with respect to different time scales and computes dimension in terms of (the logs of) length versus scale. This calculator was validated against non-fractal signals like sinusoid, box, and triangle waves, in addition to signals with known non-trivial Hausdorff dimensions and audio samples extracted from simple songs like nursery rhymes (for which we have high belief they are not fractal).

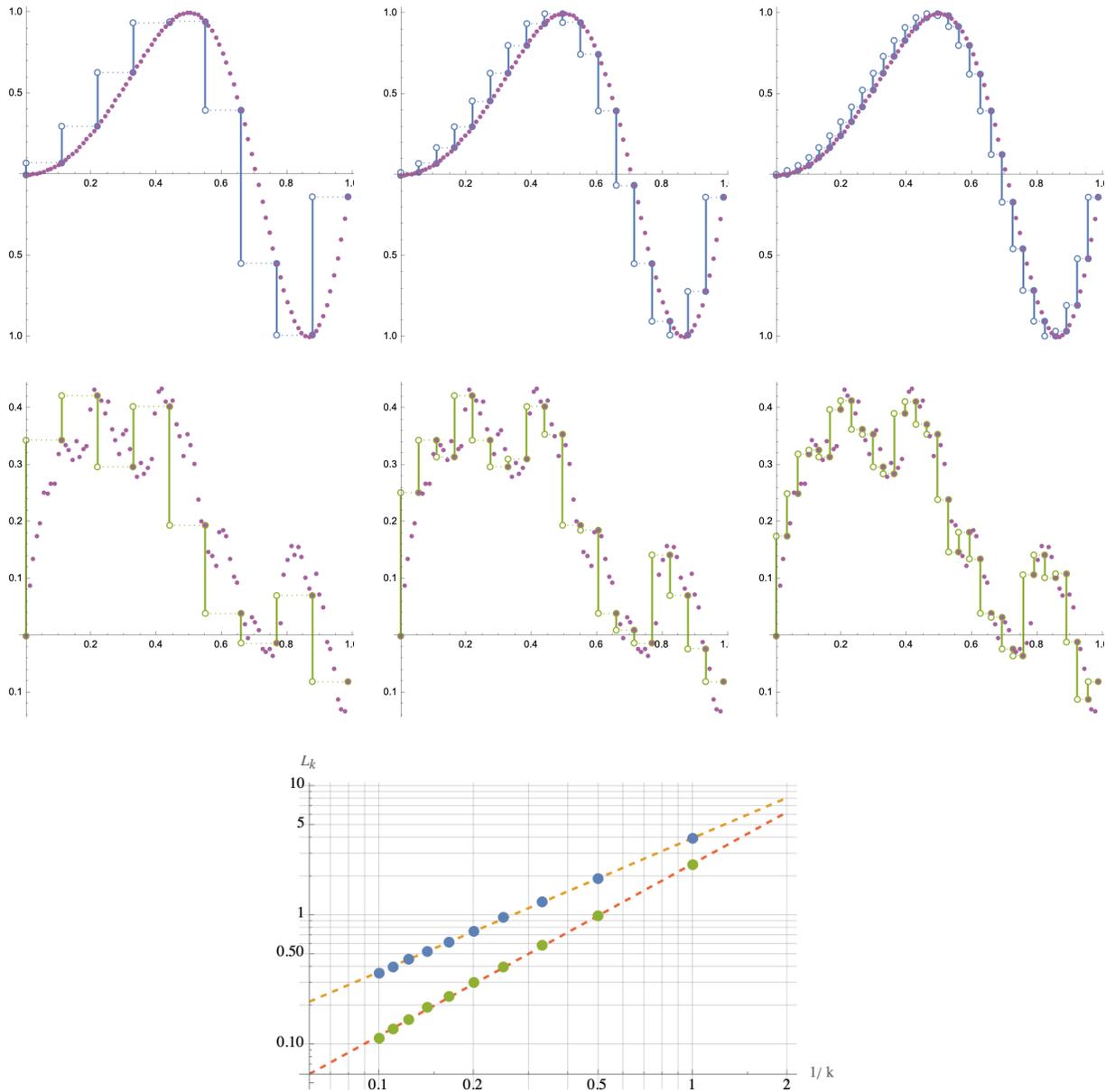

Figure 3: Computing the Hausdorff dimension of two time series by Higuchi's method. The length at a fixed scale $k$ is given in terms of absolute change in the signal over steps of size $k$, indicated by the solid vertical lines (blue and green, respectively) in the first two rows. The first signal is non-fractal, where finer measurements do not reveal additional "length," but the second Weierstrass-type signal has features that become evident as $k$ varies. As before, dimension is approximated by comparing these quantities on a log-log plot and measuring best fit slopes: the first signal yields 1.02, but the second is 1.33.

Once we were confident in the calculator, we could then test our hypothesis that the sounds that Wilson had taught us to suspect were fractal actually corresponded to fractional Hausdorff dimension. We proceeded to listen to all the songs Wilson had originally shared with Jesse in 2012 and indicated where we thought "fractal sounds" were. We then computed the Hausdorff dimensions of each of the audio files, and compared our qualitative and quantitative results:

| Song | Expect Fractal? | Fractal Dimension |
|---|---|---|
| Chant Avec Hochet (Tcoz'ungo Tzisi) | No | 1.23 |
| Le Porc - Epic | Maybe | 1.05 |
| Le Renard Aux Grandes Oreilles | Yes | 1.03 |
| Les Cornes De Rhinoceros (Hai Hu Tzi) | No | 1.03 |
| Le Melon (Tamah Tzi) | No | 1.02 |
| La Corde (Nharu Tzi) | Yes | 1.02 |
| DJamil | Yes | 1.13 |
| Xaan - Njogoy | Yes | 1.10 |
| Gainde | Yes | 1.09 |
| iRobhane | Maybe | 1.06 |
| Ngqishani, Ndabaleka | Yes | 1.02 |
| Nontyolo | No | 1.02 |
| Wavel' uSontonjane | Maybe | 1.02 |
| Medley of Baptist Sankies | No | 1.24 |
| Doption Mix | Maybe | 1.13 |
| Gone to Nineveh | Maybe | 1.04 |

In our validation process, our calculator confirmed *a priori* expectations up to two decimal places. In the data we collected, we considered a Hausdorff dimension reading of 1.02 or below to be least fractal, 1.03 to 1.07 as moderately fractal, and above 1.09 to be highly fractal.

     With a few exceptions, we found good agreement between the quantitative and qualitative: samples we expected to be fractal tended to have Hausdorff dimension higher than those for which we did not expect fractal structure. After this initial analysis, we reached out to Wilson to request a larger collection of music samples with the goal of better understanding the points of mismatch in our qualitative and quantitative analyses (for instance, in our expectations versus measurements for "Medley of Baptist Sankies"), and also to begin investigating to what extent Hausdorff dimensions might reflect particular traditions. Wilson provided us with 390 music samples over 28 albums, drawn eclectically from his multi-decade engagement with African performance cultures. As we expected from our qualitative analysis, the music samples with complex drumming patterns from Senegal produced a high Hausdorff dimension reading, as did similarly rich Congolese percussion. However, contrary to our expectations, many music samples that seemed to our ears like modern Western hip-hop (for instance, Aida Samb) which we expected to be non-fractal, produced a surprisingly high dimension reading (in the range of 1.4–1.5).

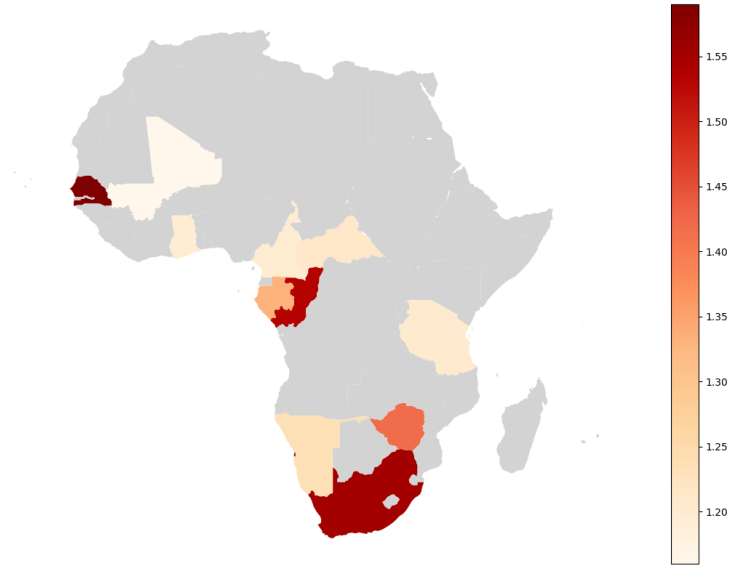

Figure 4: The maximum fractal dimension readings from our samples, according to national origin. Note that gray indicates countries with no songs in our sample. We acknowledge the limitations of this representation and look forward to its refining and problematizing in future research.

We hypothesized that modern sound-engineering and mixing techniques could account for these unexpectedly high readings, which led us to analyze a sample of contemporary Western pop music (50 songs drawn from the Fall 2022 Billboard Top 100). While we found many high measurements, we were not able to fully replicate the results from the Senegalese hip-hop sources. We also analyzed a collection of explicitly Africanist house songs from the 2000s London and New York House music scenes, which produced ratings closer to the anomalous African samples and gave us a basis to propose specific features that might account for these high Hausdorff dimensions. The initial Africanist samples can be distinguished from "modern" songs by the perceivable production techniques used in the latter---the use of synthesizers and mixing/mastering tools, like compression---and we suspect that these features might contribute to nontrivial Hausdorff dimension. We have not yet been able to verify or disprove these conjectures, but we are fascinated to continue to investigate where these high-dimensional readings come from in modern pop music.

# Looking back, looking ahead

There is an extensive literature from the last five decades on fractals and music, beginning in the 1970s with Voss and Clark's study of *1/f* scaling laws in music and sound [15, 34, 35], joining discussions of digital music [9, 10, 31], and then in work to the present [4, 6, 8, 18–21, 26, 27, 30, 32, 33, 38] as well as numerous popular accounts [5, 7, 29]. With few exceptions[26], this work is entirely focused on Western music samples, and predominantly Western scored compositional music. The majority of the analyses focus on representations of music, whether as scores [18–20,38] or via other representations of musical features [4,32,33], rather than direct analysis of waveforms or other auditory data directly available to a person hearing the music. This lack of engagement with direct audio data has been justified by claims that timbre or other features of sound might cloud or distract from "essential features" of the music [27]. Scored music—and more generally a practice of formal composition distinct from the actual creation of sound—is a rarity among the world's music traditions and practices. While the focus of prior work on primarily Western sources might reflect the priorities and interests of the authors, it might also be a natural consequence of focusing the analysis on notes on paper instead of on actual sound. The questions and methods driving productions of knowledge are bound up in the assumptions we make about the world and issues of power, status, and history.

      Our investigation verifies Wilson's hypothesis that use of fractal symmetry is not limited to African *material* cultures. How, and to what extent, the usage of fractals in Africanist music and movement cultures is distinctive among world cultures is a rich subject to explore, and not one we can do justice to on our own. We are not anthropologists or musicologists or historians, but mathematicians who are working with choreographers and performers, and who are moved by these experiential ways of interrogating mathematics. This work is part of an ongoing conversation with a wide array of people across many fields of study. We aim to continue this exchange of knowledge across cultures and disciplines on our campuses, in rehearsal studios, and beyond.

Acknowledgements

This work grew out of long-running investigations and hypotheses of Reggie Wilson, and would not have been possible without him. We thank him for extensive helpful conversations, music samples, framing, encouragement and support.  We thank company members of the *Fist and Heel Performance Group*, especially Annie Wang, and also those involved in *Moses(es)*, including Rhetta Aleong, Dwayne Brown, Yeman Brown, Paul Hamilton, Lawrence Harding, Raja Feather Kelly, Clement Mensah, and Anna Schon.  We thank Ali Tariq Younis for his assistance with developing our methodological approach and lab notebook practices.  We thank Charlotte Griffin, Susan Manning, Dawn Norfleet, Kris Petersen, Zachary Price, Tara Rodman, Jeff Streets, Deb Thomas, and S. Ama Wray for helpful discussions. This research was supported in part through funding from the University of California, Irvine's (UCI) Undergraduate Research Opportunities Program (UROP), UCI Illuminations, and National Science Foundation Grant No. DMS-1944862.

**cgonzales@carleton.edu, Dept. of Mathematics and Statistics, Carleton College.** Claudio Gómez-Gonzáles is an Assistant Professor of Mathematics at Carleton College and editorial board member of MAA FOCUS.

**svraman@uci.edu, Dept. of Mathematics, UC Irvine.** Sidhanth Raman is a fifth year Math PhD Candidate and ARCS Scholar at University of California-Irvine advised by Jesse Wolfson.

**siddharth.viswanath@yale.edu, Dept. of Computer Science, Yale University.** Siddharth Viswanath is a second year Computer Science PhD student at Yale University.

**wolfson@uci.edu, Dept. of Mathematics, UC Irvine.** Jesse Wolfson is an Associate Professor of Mathematics at the University of California-Irvine and a member of the board of directors of *Fist and Heel Performance Group*.

More detailed information regarding songs can be found here:

svraman.com/fractals-project

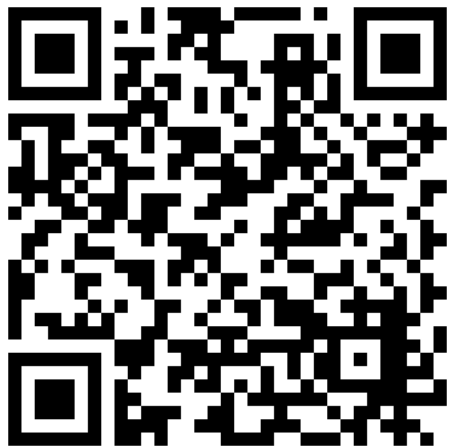